\documentclass[twocolumn,final,fleqn,envcountsect,runningheads]{my-svjour2}

\smartqed

\usepackage[T1]{fontenc}
\usepackage{textcomp}
\usepackage{mathptmx}

\usepackage{amssymb}
\usepackage{amsmath}

\journalname{The original publication is available at
www.springerlink.com}

\begin{document}

\title{Fractional conservation laws in optimal control theory\thanks{Partially
presented at FDA'06 -- \emph{2nd IFAC Workshop on Fractional
Differentiation and its Applications}, 19-21 July 2006, Porto,
Portugal (see arXiv:math.OC/0603598).}}

\author{Gast\~{a}o S. F. Frederico \and Delfim F. M. Torres}

\institute{Gast\~{a}o S. F. Frederico \at
Department of Science and Technology, University of Cape Verde\\
Praia, Santiago, Cape Verde
\and Delfim F. M. Torres (Corresponding Author: \texttt{delfim@ua.pt}) \at
Department of Mathematics, University of Aveiro\\
3810-193 Aveiro, Portugal
}

\date{}

\maketitle

\begin{abstract}
Using the recent formulation of Noether's theorem for the problems
of the calculus of variations with fractional derivatives, the
Lagrange multiplier technique, and the fractional Euler-Lagrange
equations, we prove a Noether-like theorem to the more
general context of the fractional optimal control. As a
corollary, it follows that in the fractional case the
autonomous Hamiltonian does not define anymore a conservation law.
Instead, it is proved that the fractional conservation law adds to
the Hamiltonian a new term which depends on the fractional-order
of differentiation, the generalized momentum, and the fractional
derivative of the state variable.

\keywords{Fractional derivatives \and Optimal control \and
Noether's theorem \and Conservation laws \and Symmetry}

\subclass{49K05, 26A33, 70H33.}
\end{abstract}


\section{Introduction}

The concept of symmetry plays an important role both in Physics and
Mathematics. Symmetries are described by transformations of the
system, which result in the same object after the transformation
is carried out. They are described mathematically by parameter
groups of transformations. Their importance ranges from fundamental
and theoretical aspects to concrete applications, having profound
implications in the dynamical behavior of the systems, and in
their basic qualitative properties.

Another fundamental notion in Physics and Mathematics is the one
of conservation law. Typical application of conservation laws in
the calculus of variations and optimal control is to reduce the
number of degrees of freedom, and thus reducing the problems to a
lower dimension, facilitating the integration of the differential
equations given by the necessary optimality conditions.

Emmy Noether was the first who proved, in 1918, that the notions
of symmetry and conservation law are connected: when a system
exhibits a symmetry, then a conservation law can be obtained. One
of the most important and well known illustration of this deep and
rich relation, is given by the conservation of energy in
Mechanics: the autonomous Lagrangian $L(q,\dot{q})$, correspondent
to a mechanical system of conservative points, is invariant under
time-translations (time-homogeneity symmetry),
and\footnote{Following \cite{tncdf}, we use the notation
$\partial_i L$ to denote the partial derivative of function $L$
with respect to its $i$-th argument.}
\begin{equation}
\label{eq:R}
\begin{gathered}
\frac{d}{dt} \left[ L(q,\dot{q}) - \partial_{2}
L(q,\dot{q})\cdot\dot{q} \right] = 0
\end{gathered}
\end{equation}
follows from Noether's theorem, \textrm{i.e.}, the total energy of
a conservative closed system always remain constant in time, ``it
cannot be created or destroyed, but only transferred from one form
into another''. Expression \eqref{eq:R} is valid along all the
Euler-Lagrange extremals $q(\cdot)$ of an autonomous problem of
the calculus of variations. The conservation law \eqref{eq:R} is
known in the calculus of variations as the 2nd Erdmann necessary
condition; in concrete applications, it gains different
interpretations: conservation of energy in Mechanics;
income-wealth law in Economics; first law of Thermodynamics; etc.
The literature on Noether's theorem is vast, and many extensions
of the classical results of Emmy Noether are now available for the
more general setting of optimal control (see
\cite{comIlonaIJAMAS,delfimEJC,delfimPortMath04} and references
therein). Here we remark that in all those results conservation
laws always refer to problems with integer derivatives.

Nowadays fractional differentiation plays an important role in
various fields: physics (classic and quantum mechanics,
thermodynamics, etc), chemistry, biology, economics, engineering,
signal and image processing, and control theory
\cite{Agrawal:2004b,CD:Hilfer:2000,CD:Klimek:2002}. Its origin
goes back three centuries, when in 1695 L'Hopital and Leibniz
exchanged some letters about the mathematical meaning of
$\frac{d^{n}y}{dx^{n}}$ for $n=\frac{1}{2}$. After that, many
famous mathematicians, like J.~Fourier, N.~H.~Abel, J.~Liouville,
B.~Riemann, among others, contributed to the development of the
Fractional Calculus
\cite{CD:Hilfer:2000,CD:MilRos:1993,CD:SaKiMa:1993}.

The study of fractional problems of the Calculus of Variations and
respective Euler-Lagrange type equations is a subject of current
strong research. F.~Riewe \cite{CD:Riewe:1996,CD:Riewe:1997}
obtained a version of the Euler-Lagrange equations for problems of
the Calculus of Variations with fractional derivatives, that
combines the conservative and non-conservative cases. In 2002
O.~Agrawal proved a formulation for variational problems with
right and left fractional derivatives in the Riemann-Liouville
sense \cite{CD:Agrawal:2002}. Then, these Euler-Lagrange equations
were used by D. Baleanu and T. Avkar to investigate problems with
Lagrangians which are linear on the velocities
\cite{CD:BalAv:2004}. In \cite{Klimek2001,MR1966935} fractional
problems of the calculus of variations with symmetric fractional
derivatives are considered and correspondent Euler-Lagrange
equations obtained, using both Lagrangian and Hamiltonian
formalisms. In all the above mentioned studies, Euler-Lagrange
equations depend on left and right fractional derivatives, even
when the problem depend only on one type of them. In
\cite{Klimek2005} problems depending on symmetric derivatives are
considered for which Euler-Lagrange equations include only the
derivatives that appear in the formulation of the problem. In
\cite{El-Nabulsi2005b,El-Nabulsi2005a} Riemann-Liouville
fractional integral functionals, depending on a parameter $\alpha$
but not on fractional-order derivatives of order $\alpha$, are
introduced and respective fractional Euler-Lagrange type equations
obtained. More recently, the authors have used the results of
\cite{CD:Agrawal:2002} to generalize the classical Noether's
theorem for the context of the Fractional Calculus of Variations
\cite{tncdf}. Differently from \cite{tncdf}, where the Lagrangian
point of view is considered, here we adopt an Hamiltonian point of
view. Fractional Hamiltonian dynamics is a very recent subject but
the list of publications has become already a long one due to many
applications in mechanics and physics \cite{MR2282282,MR2268869,%
El-Nabulsi1,MR2169356,MR2239336,MR2279972,MR2238509}. We extend
the previous optimal control Noether results of
\cite{delfimEJC,delfimPortMath04} to the wider context of
fractional optimal control (Theorem~\ref{thm:mainResult:FDA06}).
This is accomplished by means (i) of the fractional version of
Noether's theorem \cite{tncdf}, (ii) and the Lagrange multiplier
rule \cite{Agrawal:2004a}. As a consequence of our main result, it
follows that the ``total energy'' (the autonomous Hamiltonian) of
a fractional system is not conserved: a new expression appears
(\textrm{cf.} Corollary~\ref{cor:MainResult}) which also depends
on the fractional-order of differentiation, the adjoint variable,
and the fractional derivative of the state trajectory.


\section{Fractional Derivatives}
\label{sec:fdRL}

We briefly recall the definitions of right and left
Riemann-Liouville fractional derivatives, as well as their main
properties \cite{CD:Agrawal:2002,CD:MilRos:1993,CD:SaKiMa:1993}.

\begin{definition}
Let $f$ be a continuous and integrable function in the interval
$[a,b]$. For all $t \in [a,b]$, the left Riemann-Liouville
fractional derivative $_aD_t^\alpha f(t)$, and the right
Riemann-Liouville fractional derivative $_tD_b^\alpha f(t)$, of
order $\alpha$, are defined in the following way:
\begin{equation}
\label{eq:DFRLE}
_aD_t^\alpha f(t) =
\frac{1}{\Gamma(n-\alpha)}\left(\frac{d}{dt}\right)^{n}
\int_a^t (t-\theta)^{n-\alpha-1}f(\theta)d\theta \, ,
\end{equation}
\begin{equation}
\label{eq:DFRLD}
_tD_b^\alpha f(t) =
\frac{1}{\Gamma(n-\alpha)}\left(-\frac{d}{dt}\right)^{n}
\int_t^b(\theta - t)^{n-\alpha-1}f(\theta)d\theta \, ,
\end{equation}
where $n \in \mathbb{N}$, $n-1 \leq \alpha < n$, and $\Gamma$ is
the Euler gamma function.
\end{definition}

\begin{remark}
If $\alpha$ is an integer, then from \eqref{eq:DFRLE}
and \eqref{eq:DFRLD} one obtains the standard derivatives, that is,
\begin{gather*}
\label{eq:DU}
_aD_t^\alpha f(t) = \left(\frac{d}{dt}\right)^\alpha f(t) \, , \\
_tD_b^\alpha f(t) = \left(-\frac{d}{dt}\right)^\alpha f(t) \, .
\end{gather*}
\end{remark}

\begin{theorem} Let $f$ and $g$ be two continuous functions
on $[a,b]$. Then, for all $t \in [a,b]$, the following properties
hold:
\begin{enumerate}
\item for $p>0$, $$_aD_t^p\left(f(t)+g(t)\right)
= {_aD_t^p}f(t)+{_aD_t^p}g(t) \, ;$$
\item for $p \geq q \geq 0$,
$$_aD_t^p\left(_aD_t^{-q} f(t)\right) = {_aD_t^{p-q}}f(t) \, ;$$
\item for $p>0$, $$_aD_t^p\left(_aD_t^{-p} f(t)\right) = f(t)$$
(fundamental property of the Riemann-Liouville fractional
derivatives).
\end{enumerate}
\end{theorem}

\begin{remark}
\label{rem:der:const:NZ}
In general, the fractional derivative of a constant is not
equal to zero.
\end{remark}

\begin{remark}
The fractional derivative of order $p>0$
of function $(t-a)^\upsilon$, $\upsilon>-1$,
is given by
\begin{equation*}
_aD_t^p(t-a)^\upsilon
= \frac{\Gamma(\upsilon+1)}{\Gamma(-p+\upsilon+1)}(t-a)^{\upsilon-p} \, .
\end{equation*}
\end{remark}

\begin{remark}
When one reads ``Riemann-Liouville fractional derivative'' in the
literature, it is usually meant (implicitly) the left
Riemann-Liouville fractional derivative. In Physics, $t$ often
denotes the time-variable, and the right Riemann-Liouville
fractional derivative of $f(t)$ is interpreted as a future state
of the process $f(t)$. For this reason, right derivatives are
usually neglected in applications: the present state of a process
does not depend on the results of the future development.
Following \cite{Agrawal:2004a}, and differently from \cite{tncdf},
in this work we focus on problems with left Riemann-Liouville
fractional derivatives only. This has the advantage of simplifying
greatly the theory developed in \cite{tncdf}, making possible the
generalization of the results to the fractional optimal control
setting.
\end{remark}

We refer the interested reader in additional background on
fractional theory, to the comprehensive book
\cite{CD:SaKiMa:1993}.


\section{Preliminaries}

In \cite{CD:Agrawal:2002} a formulation of the
Euler-Lagrange equations is given for problems of the calculus of
variations with fractional derivatives.

Let us consider the following fractional problem
of the calculus of variations: to find function
$q(\cdot)$ that minimizes the integral functional
\begin{gather}
\label{Pf}
I[q(\cdot)] = \int_a^b L\left(t,q(t),{_aD_t^\alpha} q(t)\right) dt \, ,
\end{gather}
where the Lagrangian
$L :[a,b] \times \mathbb{R}^{n} \times
\mathbb{R}^{n} \rightarrow \mathbb{R}$
is a $C^{2}$ function with respect to all its arguments, and
$0 < \alpha \leq 1$.

\begin{remark}
In the case $\alpha = 1$, problem \eqref{Pf} is reduced to
the classical problem
\begin{equation*}
I[q(\cdot)] = \int_a^b L\left(t,q(t),\dot{q}(t)\right) dt
\longrightarrow \min \, .
\end{equation*}
\end{remark}

\begin{theorem}[\textrm{cf.} \cite{CD:Agrawal:2002}]
\label{Thm:FractELeq} If $q$ is a minimizer of problem \eqref{Pf},
then it satisfies the \emph{fractional Euler-Lagrange equations}:
\begin{equation}
\label{eq:eldf}
\partial_{2} L\left(t,q,{_aD_t^\alpha q}\right)
+ {_tD_b^\alpha}\partial_{3} L\left(t,q,{_aD_t^\alpha q}\right) = 0 \, .
\end{equation}
\end{theorem}

The following definition is useful in order to introduce an
appropriate concept of \emph{fractional conservation law}. We
recall that the classical Noetherian conservation laws are always
a sum of products (as assumed in \eqref{eq:somaPrd}) and that the
fractional rule for differentiation of a product, in the sense of
Riemann-Liouville, is enough complex (see \textrm{e.g.}
\cite{CD:SaKiMa:1993}). With respect to this, our operator
$\mathcal{D}_{t}^{\alpha}\left(f,g\right)$ is useful. This
operator was introduced in \cite{tncdf} and we refer the reader to
this reference for several illustrative examples and remarks. Here
we just mention that the operator
$\mathcal{D}_{t}^{\alpha}\left(f,g\right)$ has resemblances with
the classical Poisson bracket (\textrm{cf.} Remark~\ref{rem:PP}).

\begin{definition}[\textrm{cf.} \cite{tncdf}]
Given two functions $f$ and $g$ of class $C^{1}$ in the interval
$[a,b]$, we define the following operator:
\begin{equation*}
\mathcal{D}_{t}^{\alpha}\left(f,g\right)
=  f \, {_aD_t^\alpha} g - g \, {_tD_b^\alpha} f \, ,
\quad t \in [a,b].
\end{equation*}
\end{definition}

\begin{remark}
For $\alpha = 1$, operator $\mathcal{D}_{t}^{\alpha}$ is reduced to
\begin{equation*}
\begin{split}
\mathcal{D}_{t}^{1}\left(f,g\right)
&= f \, {_aD_t^1} g - g \, {_tD_b^1 f}\\
&=  f \dot{g} + \dot{f} g = \frac{d}{dt}(f g) \, .
\end{split}
\end{equation*}
\end{remark}

\begin{remark}
The linearity of the operators $_aD_t^\alpha$ and $_tD_b^\alpha$
imply the linearity of the operator $\mathcal{D}_{t}^{\alpha}$.
\end{remark}

\begin{definition}[\textrm{cf.} \cite{tncdf}]
\label{eq:fcl}
We say that $C_{f}\left(t,q,{_aD_t^\alpha q}\right)$
is a \emph{fractional conservation law} if and only if
it is possible to write $C_{f}$ in the form of a sum of products,
\begin{equation}
\label{eq:somaPrd}
C_{f}\left(t,q,d\right)
= \sum_{i=1}^{r} C_{i}^1\left(t,q,d\right) \cdot C_{i}^2\left(t,q,d\right)
\end{equation}
for some $r \in \mathbb{N}$, and for each $i = 1,\ldots,r$
the pair $C_{i}^1$ and $C_{i}^2$ satisfy one of the following relations:
\begin{equation}
\label{eq:def:lcf1}
\mathcal{D}_{t}^{\alpha}\left(C_{i}^1\left(t,q,{_aD_t^\alpha q}\right),
C_{i}^2\left(t,q,{_aD_t^\alpha q}\right)\right) = 0
\end{equation}
or
\begin{equation}
\label{eq:def:lcf2}
\mathcal{D}_{t}^{\alpha}\left(C_{i}^2\left(t,q,{_aD_t^\alpha q}\right),
C_{i}^1\left(t,q,{_aD_t^\alpha q}\right)\right) = 0
\end{equation}
along all the fractional Euler-Lagrange extremals
(\textrm{i.e.} along all the solutions of the fractional
Euler-Lagrange equations \eqref{eq:eldf}).
We then write
$\mathcal{D}\left\{{C}_{f}\left(t,q,{_aD_t^\alpha q}\right)\right\} = 0$.
\end{definition}

\begin{remark}
\label{rem:PP} For $\alpha = 1$ \eqref{eq:def:lcf1} and
\eqref{eq:def:lcf2} coincide, and
$$\mathcal{D}\left\{C\left(t,q,{_aD_t^\alpha q}\right)\right\} = 0$$
is reduced to
\begin{equation*}
\frac{d}{dt}\left\{C\left(t,q(t),\dot{q}(t)\right)\right\} = 0
\Leftrightarrow C\left(t,q(t),\dot{q}(t)\right)
\equiv \text{constant} \, ,
\end{equation*}
which is the standard meaning of \emph{conservation law},
\textrm{i.e.} a function $C\left(t,q,\dot{q}\right)$ preserved
along all the Euler-Lagrange extremals $q(t)$, $t \in [a,b]$, of
the problem. This implies that if $\left(p(t), q(t)\right)$ is a
solution to the classical Hamilton-Jacobi equations of motion,
then $C$ defines a conservation law of the Hamiltonian equations
with Hamiltonian $H$ if $\{H,C\}=0$ or $\{C,H\}=0$, where
$\{\cdot,\cdot\}$ denotes the canonical Poisson bracket operator.
In the more general fractional context, the Hamilton-Jacobi
equations were recently derived in
\cite{CD:Jumarie:2007a,CD:Jumarie:2007b}.
\end{remark}

\begin{definition}[\textrm{cf.} \cite{tncdf}]
\label{def:invadf}
Functional \eqref{Pf} is said to be invariant
under the one-parameter group of infinitesimal transformations
\begin{equation}
\label{eq:tinf2}
\begin{cases}
\bar{t} = t + \varepsilon\tau(t,q) + o(\varepsilon) \, ,\\
\bar{q}(\bar{t}) = q(t) + \varepsilon\xi(t,q) + o(\varepsilon) \, ,\\
\end{cases}
\end{equation}
if, and only if,
\begin{multline}
\label{eq:invdf}
\int_{t_{a}}^{t_{b}} L\left(t,q(t),{_{t_a}D_t^\alpha q(t)}\right) dt \\
= \int_{\bar{t}(t_a)}^{\bar{t}(t_b)} L\left(\bar{t},\bar{q}(\bar{t}),
{_{\bar{t}_a}D_{\bar{t}}^\alpha \bar{q}(\bar{t})}\right) d\bar{t}
\end{multline}
for any subinterval $[{t_{a}},{t_{b}}] \subseteq [a,b]$.
\end{definition}

\begin{remark}
Having in mind that condition \eqref{eq:invdf}
is to be satisfied for any subinterval $[{t_{a}},{t_{b}}] \subseteq [a,b]$,
we can rid off the integral signs in \eqref{eq:invdf}.
This is done in the new Definition~\ref{def:inv:gt}.
\end{remark}

The next theorem provides an extension of the classical Noether's
theorem to Fractional Problems of the Calculus of Variations.

\begin{theorem}[\textrm{cf.} \cite{tncdf}]
\label{theo:tndf} If functional \eqref{Pf} is invariant
under \eqref{eq:tinf2}, then
\begin{multline*}
\left[L\left(t,q,{_aD_t^\alpha q}\right)
- \alpha\partial_{3} L\left(t,q,{_aD_t^\alpha q}\right)
\cdot{_aD_t^\alpha q} \right] \tau(t,q) \\
+ \partial_{3} L\left(t,q,{_aD_t^\alpha q}\right) \cdot \xi(t,q)
\end{multline*}
is a fractional conservation law (\textrm{cf.} Definition~\ref{eq:fcl}).
\end{theorem}


\section{Main Results}
\label{sec:MR}

Using Theorem~\ref{theo:tndf}, we obtain here
a  Noether's Theorem for
the fractional optimal control problems
introduced in \cite{Agrawal:2004a}:
\begin{gather}
\label{eq:JO}
I[q(\cdot),u(\cdot)] =\int_a^b L\left(t,q(t),u(t)\right) dt
\longrightarrow \min \, , \\
_aD_t^\alpha q(t)=\varphi\left(t,q(t),u(t)\right) \, , \notag
\end{gather}
together with the initial condition $q(a)=q_a$. In problem
\eqref{eq:JO}, the Lagrangian $L : [a,b]\times \mathbb{R}^{n}
\times \mathbb{R}^{m} \rightarrow \mathbb{R}$ and the velocity
vector $\varphi : [a,b] \times \mathbb{R}^{n} \times
\mathbb{R}^{m} \rightarrow\mathbb{R}^n$ are assumed to be $C^{1}$
functions with respect to all the arguments. In agreement with the
calculus of variations, we also assume that the admissible control
functions take values on an open set of $\mathbb{R}^m$.

\begin{definition}
A pair $(q(\cdot),u(\cdot))$ satisfying the fractional control
system $_aD_t^\alpha q(t)=\varphi\left(t,q(t),u(t)\right)$ of
problem \eqref{eq:JO}, $t \in [a,b]$, is called a \emph{process}.
\end{definition}

\begin{theorem}[\textrm{cf.} (13)-(15) of \cite{Agrawal:2004a}]
\label{th:P}
If $(q(\cdot),u(\cdot))$ is an optimal process for
problem \eqref{eq:JO}, then there exists a co-vector function
$p(\cdot)$ such that the following conditions hold:
\begin{itemize}
\item the Hamiltonian system
\begin{equation*}
\label{eq:Ham}
\begin{cases}
_aD_t^\alpha q(t)&=\partial_4 {\cal H}(t, q(t), u(t),p(t)) \, , \\
_tD_b^\alpha p(t) &= \partial_2{\cal H}(t,q(t),u(t), p(t)) \, ;
\end{cases}
\end{equation*}
\item the stationary condition
\begin{equation*}
\label{eq:CE}
 \partial_3 {\cal H}(t, q(t), u(t), p(t))=0 \, ;
\end{equation*}
\end{itemize}
with the Hamiltonian ${\cal H}$ defined by
\begin{equation}
\label{eq:H} {\cal H}\left(t,q,u,p\right) = L\left(t,q,u\right)
+ p \cdot \varphi\left(t,q,u\right) \, .
\end{equation}
\end{theorem}

\begin{remark}
In classical mechanics, the Lagrange multiplier $p$ is called the \emph{generalized momentum}.
In the language of optimal control, $p$ is known as the \emph{adjoint variable}.
\end{remark}

\begin{definition}
\label{def:extPont} Any triplet $(q(\cdot),u(\cdot),p(\cdot))$
satisfying the conditions of Theorem~\ref{th:P} will be called a
\emph{fractional Pontryagin extremal}.
\end{definition}

For the fractional problem of the calculus of variations \eqref{Pf}
one has $\varphi(t,q,u)=u \Rightarrow {\cal H} = L + p \cdot u$,
and we obtain from Theorem~\ref{th:P} that
\begin{gather*}
_aD_t^\alpha q = u \, ,\\
_tD_b^\alpha  p  = \partial_2 L \, ,\\
 \partial_3 {\cal H} = 0 \Leftrightarrow
 p= - \partial_3 L \Rightarrow
{_tD_b}^\alpha p= - _tD_b^\alpha \partial_3 L \, .
\end{gather*}
Comparing the two expressions for $_tD_b^\alpha p$, one arrives to the
Euler-Lagrange differential equations \eqref{eq:eldf}:
$\partial_2 L = - _tD_b^\alpha \partial_3 L$.

We define the notion of invariance for problem \eqref{eq:JO}
in terms of the Hamiltonian, by introducing the augmented functional
as in \cite{Agrawal:2004a}:
\begin{multline}
\label{eq:J} J[q(\cdot),u(\cdot),p(\cdot)] \\
= \int_a^b \left[{\cal
H}\left(t,q(t),u(t),p(t)\right)-p(t) \cdot {_aD_t}^\alpha q(t)\right]dt \, ,
\end{multline}
where ${\cal H}$ is given by \eqref{eq:H}.

\begin{remark}
Theorem~\ref{th:P} is easily obtained applying the necessary optimality
condition \eqref{eq:eldf} to problem \eqref{eq:J}.
\end{remark}

\begin{definition}
\label{def:inv:gt} A fractional optimal control problem \eqref{eq:JO} is
said to be invariant under the $\varepsilon$-parameter local group
of transformations
\begin{equation}
\label{eq:trf:inf}
\begin{cases}
\bar{t} = t+\varepsilon\tau(t, q(t), u(t), p(t)) + o(\varepsilon) \, , \\
\bar{q}(\bar{t}) = q(t)+\varepsilon\xi(t, q(t), u(t), p(t)) + o(\varepsilon) \, , \\
\bar{u}(\bar{t}) = u(t)+\varepsilon\sigma(t, q(t), u(t), p(t)) + o(\varepsilon) \, , \\
\bar{p}(\bar{t}) = p(t)+\varepsilon \zeta(t, q(t), u(t), p(t))+ o(\varepsilon) \, , \\
\end{cases}
\end{equation}
if, and only if,
\begin{multline}
\label{eq:condInv}
\left[{\cal H}(\bar{t},\bar{q}(\bar{t}),\bar{u}(\bar{t}),\bar{p}(\bar{t}))
-\bar{p}(\bar{t}) \cdot  {_{\bar{a}}D_{\bar{t}}}^\alpha \bar{q}(\bar{t})\right] d\bar{t} \\
=\left[{\cal H}(t,q(t),u(t),p(t))-p(t) \cdot {_aD_t}^\alpha q(t)\right] dt \, .
\end{multline}
\end{definition}

\begin{theorem}[Fractional Noether's theorem]
\label{thm:mainResult:FDA06}
If the fractional optimal control problem \eqref{eq:JO} is invariant
under \eqref{eq:trf:inf}, then
\begin{equation}
\label{eq:tndf:CO}
\left[ {\cal H} - \left(1 - \alpha\right) p(t)
\cdot {_aD_t}^\alpha q(t) \right] \tau - p(t) \cdot \xi
\end{equation}
is a fractional conservation law, that is,
\begin{equation*}
\mathcal{D}\left\{
\left[ {\cal H} - \left(1 - \alpha\right) p(t)
\cdot {_aD_t}^\alpha q(t) \right] \tau - p(t) \cdot \xi
\right\} = 0
\end{equation*}
along all the fractional Pontryagin extremals.
\end{theorem}

\begin{remark}
For $\alpha = 1$, the fractional optimal control problem
\eqref{eq:JO} is reduced to the classical optimal control problem
\begin{gather*}
I[q(\cdot),u(\cdot)] =\int_a^b L\left(t,q(t),u(t)\right) dt
\longrightarrow \min \, , \\
\dot{q}(t)=\varphi\left(t,q(t),u(t)\right) \, ,
\end{gather*}
and we obtain from Theorem~\ref{thm:mainResult:FDA06}
the optimal control version of Noether's theorem \cite{delfimEJC}:
invariance under a one-parameter group of transformations
\eqref{eq:trf:inf} imply that
\begin{equation}
\label{eq:H9} C(t,q,u,p)={\cal H}(t,q,u,p)\tau-p\cdot \xi
\end{equation}
is constant along any Pontryagin extremal
(one obtains \eqref{eq:H9} from \eqref{eq:tndf:CO}
setting $\alpha = 1$).
\end{remark}

\begin{proof}
The fractional conservation law \eqref{eq:tndf:CO}
is obtained applying Theorem~\ref{theo:tndf}
to the augmented functional \eqref{eq:J}. \qed
\end{proof}

Theorem~\ref{thm:mainResult:FDA06} provides a new interesting insight
for the fractional autonomous variational problems.
Let us consider the autonomous fractional optimal control problem,
\textrm{i.e.} the situation when the Lagrangian $L$
and the fractional velocity vector $\varphi$
do not depend explicitly on time $t$:
\begin{equation}
\label{eq:FOCP:CA}
\begin{gathered}
I[q(\cdot),u(\cdot)] =\int_a^b L\left(q(t),u(t)\right) dt
\longrightarrow \min \, , \\
_aD_t^\alpha q(t)=\varphi\left(q(t),u(t)\right) \, .
\end{gathered}
\end{equation}

\begin{corollary}
\label{cor:MainResult}
For the autonomous problem \eqref{eq:FOCP:CA} the following
fractional conservation law holds:
\begin{equation}
\label{eq:ConsHam:alpha}
\mathcal{D}\left\{
{\cal H} - \left(1 - \alpha\right) p(t)
\cdot {_aD_t}^\alpha q(t)
\right\} = 0 \, .
\end{equation}
\end{corollary}

\begin{remark}
In the classical framework of optimal control theory one has $\alpha = 1$
and our operator $\mathcal{D}$ coincides with $\frac{d}{dt}$.
We then get from \eqref{eq:ConsHam:alpha}
the classical result: the Hamiltonian ${\cal H}$ is a preserved
quantity along any Pontryagin extremal of the problem.
\end{remark}

\begin{proof}
The Hamiltonian ${\cal H}$ does not depend explicitly on time,
and it is easy to check that \eqref{eq:FOCP:CA}
is invariant under time-translations: invariance condition \eqref{eq:condInv}
is satisfied with $\bar{t} = t + \varepsilon$, $\bar{q}(\bar{t}) = q(t)$,
$\bar{u}(\bar{t}) = u(t)$ and $\bar{p}(\bar{t}) = p(t)$. In fact,
given that $d\bar{t} = dt$, \eqref{eq:condInv} holds trivially
proving that ${_{\bar{a}}D_{\bar{t}}}^\alpha \bar{q}(\bar{t})
= {_aD_t}^\alpha q(t)$:
\begin{equation*}
\begin{split}
_{\bar{a}} & D_{\bar{t}}^\alpha \bar{q}(\bar{t}) \\
&= \frac{1}{\Gamma(n-\alpha)}\left(\frac{d}{d\bar{t}}\right)^{n}
\int_{\bar{a}}^{\bar{t}} (\bar{t}-\theta)^{n-\alpha-1}\bar{q}(\theta)d\theta \\
&= \frac{1}{\Gamma(n-\alpha)}\left(\frac{d}{dt}\right)^{n}
\int_{a + \varepsilon}^{t+\varepsilon} (t + \varepsilon-\theta)^{n-\alpha-1}\bar{q}(\theta)d\theta \\
&= \frac{1}{\Gamma(n-\alpha)}\left(\frac{d}{dt}\right)^{n}
\int_{a}^{t} (t-s)^{n-\alpha-1}\bar{q}(s + \varepsilon)ds \\
&= {_{a}D_{t}}^\alpha \bar{q}(t + \varepsilon) = {_{a}D_{t}}^\alpha \bar{q}(\bar{t}) \\
&= {_{a}D_{t}}^\alpha q(t) \, .
\end{split}
\end{equation*}
Using the notation in \eqref{eq:trf:inf}, one has $\tau = 1$ and
$\xi = \sigma = \zeta = 0$. Conclusion \eqref{eq:ConsHam:alpha}
follows from Theorem~\ref{thm:mainResult:FDA06}. \qed
\end{proof}


\section{Illustrative Examples}

We begin by illustrating our results with two Lagrangians that do not depend explicitly on the time variable $t$. These two examples are borrowed
from \cite[\S 4.1]{Agrawal:2004a} and \cite[\S 3.1]{MR2279972},
where the authors write down the respective fractional
Euler-Lagrange equations. Here, we use our
Corollary~\ref{cor:MainResult} to obtain new fractional
conservation laws.

\begin{example}
We begin by considering a simple fractional problem of the
calculus of variations (see \cite[Example~1]{CD:Agrawal:2002} and
\cite[\S 3.1]{MR2279972}):
\begin{equation}
I[q(\cdot)] = \frac{1}{2}\int_0^1 \left(_{0}D_{1}^{\alpha}q(t)
\right)^2dt \longrightarrow \min \, , \quad \alpha > \frac{1}{2} \, .
\end{equation}
Equation \eqref{eq:H} takes the form
\begin{equation}
\label{eq:cons:Energ:Ex1} {\cal H}=-\frac{1}{2}p^2\,.
\end{equation}
We conclude from Corollary~\ref{cor:MainResult} that
\begin{equation}
\label{eq:lc:ex1}
\frac{p^2}{2}(1-2\alpha)
\end{equation}
is a fractional conservation law.
\end{example}

\begin{example}
Let us now consider the following fractional optimal control
problem \cite[\S 4.1]{Agrawal:2004a}:
\begin{gather}
\label{Po} I[q(\cdot)] = \frac{1}{2}\int_0^1 \left[q^2(t)+
u^2(t)\right] dt
\longrightarrow \min \, ,\\
_{0}D_{1}^{\alpha}q(t)=-q(t)+u(t) \, , \notag
\end{gather}
under the initial condition $q(0)=1$. The Hamiltonian ${\cal H}$
\eqref{eq:H} has the form
$${\cal H}=\frac{1}{2}\left(q^2+
u^2\right)+p(-q+u).
$$
From Corollary~\ref{cor:MainResult} it follows that
\begin{equation}
\label{eq:cons:Energ:Ex2} \frac{1}{2}\left(q^2+u^2\right)+\alpha
p(-q+u)
\end{equation}
is a fractional conservation law.
\end{example}

For $\alpha = 1$, the fractional conservation laws
\eqref{eq:lc:ex1} and \eqref{eq:cons:Energ:Ex2} give conservation of energy.

Finally, we give an example of an optimal control problem with three state variables and two controls ($n=3$, $m=2$). The problem is inspired in \cite[Example~2]{comPauloLituania05}.

\begin{example}
We consider the following fractional optimal control problem:
\begin{eqnarray}
\int_a^b \left(u_1(t)^2+u_2(t)^2\right)\mathrm{d}t \longrightarrow \min, \label{ex:kc}\\
\left\{\begin{array}{l}
_aD_t^\alpha q_1(t)=u_1(t) \cos (q_3(t)),\\
_aD_t^\alpha q_2(t)=u_1(t) \sin (q_3(t)),\\
_aD_t^\alpha q_3(t)=u_2(t) .
\end{array}\right. \label{ex:kc:2}
\end{eqnarray}
For $\alpha = 1$ the control system \eqref{ex:kc:2} serves as model for the kinematics of a car and \eqref{ex:kc}-\eqref{ex:kc:2} reduces to
Example~2 of \cite{comPauloLituania05}.
From Corollary~\ref{cor:MainResult} one gets that
\begin{multline*}
u_1^2 + u_2^2 + p_1 \left(u_1 \cos(q_3) - (1-\alpha) _aD_t^\alpha q_1\right) \\
+ p_2 \left(u_1 \sin(q_3) - (1-\alpha) _aD_t^\alpha q_2\right)
+ p_3 \left(u_2 - (1-\alpha) _aD_t^\alpha q_3\right)
\end{multline*}
is a fractional conservation law.
\end{example}

Main difficulty of our approach is related with the computation of the invariance transformations. To illustrate this issue,
let us consider problem \eqref{eq:JO} with
$$
L(t,q,u)=L(t,u) \, , \quad \varphi(t,q,u)=\varphi(t,u) \, .
$$
In the classical case,
since $q$ does not appear both in $L$ and $\varphi$,
such a problem is trivially invariant under translations on the variable $q$, \textrm{i.e.} condition \eqref{eq:condInv} is verified for $\alpha = 1$ with $\bar{t}=t$,
$\bar{q}(t)=q(t)+\varepsilon$, $\bar{u}(\bar{t})=u(t)$ and
$\bar{p}(\bar{t})=p(t)$. In the fractional case this is not in general true: we have $d\bar{t}=dt$, but condition \eqref{eq:condInv} is not satisfied since
${_{\bar{a}}D^{\alpha}_{\bar{t}}}\bar{q}(\bar{t})
={_aD^{\alpha}_{t}}q(t) + {_aD^\alpha_{t}} \epsilon$
and the second term on the right-hand side is in general not equal to zero (Remark~\ref{rem:der:const:NZ}).


\section{Conclusions}
\label{sec:Conc}

The fractional Euler-Lagrange equations are a subject of strong
current study
\cite{CD:Agrawal:2002,CD:BalAv:2004,El-Nabulsi2005b,%
El-Nabulsi2005a,Klimek2001,MR1966935,Klimek2005,CD:Riewe:1996,CD:Riewe:1997}
because of its numerous applications. In \cite{tncdf} a fractional
Noether's theorem is proved.

The fractional Hamiltonian perspective is a quite recent subject,
being investigated in a serious of publications
\cite{MR2282282,MR2268869,El-Nabulsi1,Muslih:2005,%
MR2169356,MR2239336,MR2279972,MR2238509}. One can say, however,
that the fractional variational theory is still in its childhood.
Much remains to be done. This is particularly true in the area of
fractional optimal control where results are a rarity. The main
study of fractional optimal control problems seems to be
\cite{Agrawal:2004a}, where the Euler-Lagrange equations for
fractional optimal control problems (Theorem~\ref{th:P}) are
obtained, using the traditional approach of the Lagrange
multiplier rule. Here we use the Lagrange multiplier technique to
derive, from the results in \cite{tncdf}, a new Noether-type
theorem for fractional optimal control systems. Main result
generalizes the results of \cite{delfimEJC}. As an application, we
have considered the fractional autonomous problem, proving that
the Hamiltonian defines a conservation law only in the integer
case $\alpha = 1$.


\begin{acknowledgements}
GF was supported by IPAD (Instituto Portugu\^{e}s de Apoio ao
Desenvolvimento); DT by CEOC (Centre for Research on Optimization
and Control) through FCT (Portuguese Foundation for Science and
Technology), cofinanced by the European Community fund
FEDER/POCI 2010. We are grateful to Professor Tenreiro Machado for
drawing our attention to the \emph{2nd IFAC Workshop on Fractional
Differentiation and its Applications}, 19--21 July, 2006, Porto,
Portugal, and for encouraging us to write the present work.
Inspiring discussions with Jacky Cresson at the
Universit\'{e} de Pau et des Pays de l'Adour, France,
are very acknowledged.
\end{acknowledgements}



\end{document}